\documentclass[onefignum,onetabnum,oneeqnum,onealgnum,reqno]{siamart251216}
\usepackage{amsmath,amsfonts}
\usepackage{amssymb}
\usepackage{color}
\usepackage{float}
\usepackage{isomath}
\usepackage{graphicx}
\usepackage{hyperref}
\usepackage{orcidlink}

\usepackage{tikz}
\usetikzlibrary{math}
\usepackage[labelfont=bf,textfont=it,font=small]{caption}
\usepackage{subcaption}

\usepackage{algorithm}
\usepackage{algpseudocode} 

\definecolor{darkgray}{rgb}{0.8, 0.8, 0.8}

\algrenewcommand{\algorithmiccomment}[2][.6\linewidth]{\leavevmode\hfill\makebox[#1][l]{$\triangleright$ #2}}
\algnewcommand{\algorithmicgoto}{\textbf{goto}}
\algnewcommand{\Goto}[1]{\algorithmicgoto~\ref{#1}}%
\algnewcommand\algorithmicvariables{\textbf{Local variables:}}%
\algnewcommand\Variables{\item[\algorithmicvariables]}%

\algblock{ParFor}{EndParFor}
\algnewcommand\algorithmicparfor{\textbf{parallel for}}
\algnewcommand\algorithmicpardo{\textbf{do}}
\algnewcommand\algorithmicendparfor{\textbf{end\ parallel for}}
\algrenewtext{ParFor}[1]{\algorithmicparfor\ #1\ \algorithmicpardo}
\algrenewtext{EndParFor}{\algorithmicendparfor}

\algnewcommand{\LineComment}[1]{\vspace{3pt}\Statex \(\triangleright\)\(\triangleright\) \textcolor{blue}{#1} \vspace{3pt}}
\newcommand{\COMMENT}[2][.47\linewidth]{%
    \leavevmode\hfill\makebox[#1][l]{$\triangleright$~#2}}

\renewcommand{\Re}{\mathbb{R}}

\newcommand{\R}{\mathbb{R}} 

\newtheorem{rem}{Remark}

\renewcommand{\vec}{\vectorsym}
\newcommand{\mat}{\matrixsym}

\renewcommand{\eqref}[1]{Equation (\ref{#1})}

\newcommand{\algG}{\texttt{G}}
\newcommand{\algTG}{\texttt{TG}}
\newcommand{\algMV}{\texttt{MV}}
\newcommand{\algMW}{\texttt{MW}}
\newcommand{\algIMG}{\texttt{IMG}}
\newcommand{\algMG}{\texttt{MG}}
\newcommand{\algFMG}{\texttt{FMG}}

\newcommand{\M}{\mathcal{M}}
\newcommand{\F}{E}

\title{Multigrid Primer: Basic Principles\vspace{1em}}

\author{\orcidlinki{Stephen F. McCormick}{0000-0002-6327-790X}\thanks{University of Colorado, Boulder (\email{stephen.mccormick@colorado.edu})}
  \and \orcidlinki{Rasmus Tamstorf}{0000-0001-7649-8587}\thanks{Independent Researcher (\email{rt@acm.org})}}
\date{}

\begin{document}
\maketitle

\begin{abstract}
The goal of this primer is to provide a relatively short exposition of the basics of multigrid methods, simplified by focusing on fundamental concepts in a {\em variational} setting. This is done by way of a quadratic {\em energy} minimization formulation of symmetric positive definite linear equations arising from the discretization of elliptic partial differential equations. This focus provides an alternate viewpoint to other expositions and, more importantly, it enables a simplification of the development while clarifying the principles that lead to effective algorithms. The development begins with multigrid as an {\em iterative solver} exemplified by the so-called V-cycle. It then introduces the full multigrid method as a {\em direct solver} in the sense that it is aimed directly at the source of the matrix equations. In this way, full multigrid attempts to achieve discretization-level accuracy at a cost comparable to that of a few matrix multiplies on the finest level.
\end{abstract}



\section{Preamble}
\label{sec:pre}

Multigrid (MG) methods involve a hierarchy of discrete levels used to efficiently compute an accurate approximation of an equation on the finest level. It may seem natural to start by computing an approximation of a cruder form of the equation on the coarsest level, using the result as an initial guess for the iteration on the next level, and so on until the finest level is reached and processed. Such a {\em nested iteration} approach makes sense because the relatively inexpensive computations on the coarse levels might hopefully supply a good way to start iterating on the finest level.

On the other hand, it may seem strange to think that it is necessary to return to the coarse levels for further processing. Yet, we aim to show here that this is really quite natural. In fact, in a broad sense, we humans operate effectively in such a multilevel way in our everyday lives. For example, suppose that someone just walked in to join a meeting in a large room with many seated participants and many other empty chairs. It is normal to go back and forth looking at individual seats and taking broader views of the layout before selecting one's own chair. Think also of an artist or writer who works closely on detail but often backs up to the overall view to see where local changes have triggered necessary adjustments to the bigger picture, only then to returns to adjust the local details that the larger changes have triggered. 

So it goes with multigrid. We show why such a back-and-forth process across the hierarchy is used by a multigrid solver and why it ensures efficiency.

\section{Introduction}
\label{sec:intro}

The purpose here is to provide a brief introduction to multigrid methods together with some practical principles and heuristics that can be used to guide their design. Some of this presentation has been taken from the supplemental material available with \cite{Tamstorf2015}, while most of the following concepts can be studied further in \cite{briggs00}, \cite{haase-langer}, \cite{mccormick84}, \cite{trottenberg00}, and \cite{xu2017}. Our focus in what follows is on symmetric positive definite (SPD) matrix equations that arise from the discretization of elliptic partial differential equations (PDEs). Although multigrid and, more generally, multilevel methods have much broader applicability, this case allows for a much simpler development and clarification of the basic principles underlying multigrid methodology. 

This Primer focuses on general principles rather than detailed prescriptions for specific problems, emphasizing how understanding these principles can guide the development of effective multigrid methods for diverse applications.

For a glimpse of the simplicity of the variational form of multigrid methods, note first that they have two essential components: {\em relaxation} (also called {\em smoothing}) and {\em coarsening} (often called {\em coarse-grid correction}). Typically, relaxation is a fairly simple iterative method with the responsibility of removing oscillatory errors in the solution, thereby making the error {\em smooth} in a sense that we describe later, while coarsening the equations to cruder levels effectively provides a way to deal with these smooth errors. More generally, coarsening is there to deal with the broader errors in the current approximation that causes relaxation to stall. For multigrid to work optimally, these two components must be carefully designed to complement each other. An underlying theme of this Primer is the focus on design conditions that guarantee such complementarity.

In the variational setting, the goal of multigrid is effectively to minimize a functional $F(u)$ that maps some space $U$ of vector-valued functions to the real line $\Re$. Speaking formally, we might consider a gradient method that seeks to replace a current approximation $u$ by $u - s \nabla F(u)$, where $s$ is a real number chosen in an attempt to minimize $F(u - s \nabla F(u))$. As we show below, this is the form that relaxation takes (except that we take a simpler choice for $s$). Suppose that we then want to update the new $u$ by finding a $v$ in a vector-valued function space $V$ that attempts to improve $u$ by reducing its functional value. If $P$ is a mapping from $V$ into $U$, then the task might be to replace $u$ by $u - Pv$, where $v \in V$ is chosen in an attempt to minimize $F(u - Pv)$. This is the form that coarsening takes. Put together, a multigrid cycle with two grids is therefore described {\em formally} as follows:
\vspace{0.5em}
\begin{itemize}
\item {\em Relaxation}. Replace $u$ by $u - s \nabla F(u)$, where $s \in \Re$ minimizes $F(u - s \nabla F(u))$.
\item {\em Coarsening}. Replace $u$ by $u - Pv$, where $v \in V$ minimizes $F(u - Pv)$.
\end{itemize}
\vspace{0.5em}
This hopefully shows how simple multigrid is when viewed in the abstract. But, of course, it raises many question, not the least of which is how these minimizations are supposed to be computed. We address these questions in some detail here, but the point is that the basic form of multigrid is fairly simple.

To better understand the interplay between relaxation and coarsening, we need to spend a fair bit of time developing notation, basic concepts, and other preliminaries before attempting a simple description of a multigrid method. We start in section~\ref{sec:preliminaries} with the assumptions on the hierarchy of discrete matrix equations that are made available to the multigrid algorithms. We continue in section~\ref{sec:smoothing} by characterizing the relaxation method that we have in mind together with their basic properties. After section~\ref{sec:coarseCorrection} that introduces coarse-grid correction, section~\ref{sec:2grid}  develops a two-grid cycle, thus leading into the construction of so-called V-cycles and W-cycles in section~\ref{sec:solvers}. These cycles are used in section~\ref{sec:ims} to construct {\em iterative} multigrid methods with
uniformly bounded convergence factors {\em per cycle} at a theoretical cost proportional to the number of degrees of freedom on the finest level. Section~\ref{sec:FMG} introduces full multigrid (FMG) viewed as a {\em direct} method in the sense that its target is the PDE: it uses a special fixed cycle that attempts to achieve discretization-level accuracy at the truly optimal {\em total} cost equivalent to that of just a few finest-level relaxation sweeps. We end with a short summary in section~\ref{sec:summary}.

\section{Preliminaries}
\label{sec:preliminaries}
We use bold font for matrices and vectors, typewriter font for algorithms, and italic font for functions and scalars. For example, $\mat{G}$ signifies the error propagation matrix for a generic stationary linear relaxation method, $\algG^h\left({\vec v}^h, {\vec g}^h,\mu,h\right)$ denotes the relaxation algorithm, and G is shorthand for relaxation. We write the PDE as $Lu=f$, where $L$ denotes the operator, $u$ the solution, and $f$ the source term. For consistency, our matrix equations appear as expressions like $\mat{L}^h\vec{u}^h=\vec{f}^h$, where the superscripts specify the grid level. Note that we use this notation to be consistent with the PDE rather than the more traditional notation $\mat{A}\vec{x}=\vec{b}$ for a linear system. 

We start with some simplifying assumptions that hold in all that follows.\footnote{The assumptions we make are only to simplify the discussion. Multigrid and, more generally, multilevel methods (that might not even be based on grids) have much broader applicability.} Assume that a hierarchy of SPD matrices ${\mat L}^h \in \Re^{n \times n}$ and intergrid transfer operators called {\em interpolation} or {\em prolongation} ${\mat P}_{2h}^h$ and {\em restriction} ${\mat R}_h^{2h}$ have been determined on a sequence of coarse-to-fine grid levels that reflect increasingly accurate discretizations of a PDE. The assumption about these discretizations is that they are based on the simple case of a finest-level logically rectangular grid in two or three dimensions with uniform mesh size $h$ in all coordinate directions. A further simplification is that the coarse levels are constructed by eliminating every other grid line in each coordinate direction. Accordingly, we use $2h$ to denote the mesh size of a coarsening of a given grid $h$, where $H$ is used for the coarsest grid in the hierarchy. We do not explicitly specify here how these discretizations are to be constructed so that the discussion applies to virtually any procedure. However, it may be best to think in terms of the Rayleigh-Ritz {\em variational} approach (cf. \cite{strang-fix}) based on lowest-order finite elements. In what follows, these constructs on all levels of the hierarchy are assumed to be internally available to the algorithms that we introduce. We refer to them as being ``Given''. 

\section{Relaxation}
\label{sec:smoothing}
We begin the development of multigrid by studying relaxation.
The target problem is to find ${\vec u}^h \in\R^n$ such that
\begin{equation}
{\mat L}^h{\vec u}^h = {\vec f}^h,
\label{equation}
\end{equation}
where ${\vec f}^h \in\R^n$ is a given source term. The coarse-grid correction form of multigrid described later involves an equation on grid $2h$ where the source term is different than than the ${\vec f}^{2h}$ that comes directly from the discretization. To accommodate this change, we view the ${\vec u}^h$ and ${\vec f}^h$ in (\ref{equation}) as given and fixed, and introduce the following related equation where the ${\vec v}^h$ and ${\vec g}^h$ are allowed to change in the algorithm:
\begin{equation}
{\mat L}^h{\vec v}^h = {\vec g}^h.
\label{resid-eq}
\end{equation}
While (\ref{equation}) is our ultimate target, we focus first on (\ref{resid-eq}) in our discussions of relaxation and basic multigrid. This is of course just a change in notation, but it is useful because it allows us to emphasize the difference between treating (\ref{equation}) directly and treating it indirectly by way of (\ref{resid-eq}). This indirect approach involves setting up (\ref{resid-eq}) so that its solution is a correction to the current approximate solution of (\ref{equation}). We return to (\ref{equation}) when we introduce FMG that treats it directly. We first discuss our choice of relaxation in some depth.

\subsection{Richardson iteration} The matrix equation in (\ref{resid-eq}) can be solved using any of a number of simple iterative methods like Jacobi or Gauss-Seidel relaxation~(cf. \cite[Section 10.1]{Golub83}). However, in all that follows,
we restrict our attention to Richardson's iteration given by
\begin{equation}
{\vec v}^h \leftarrow {\vec v}^h - \frac{1}{\|{\mat L}^h\|}\left({\mat L}^h{\vec v}^h - {\vec g}^h\right),
\label{error}
\end{equation}
where $\|\cdot\|$ denotes the matrix norm induced by the Euclidean vector norm $\|\cdot\|$. (By the notation in (\ref{error}), we mean that the current approximation, ${\vec v}^h$, to the solution of (\ref{resid-eq}) is replaced by the expression involving the old ${\vec v}^h$ to the right of the arrow.) As with many other methods, Richardson relaxation corrects the current approximation by a term involving the {\em residual} ${\vec r}^h \equiv {\mat L}^h{\vec v}^h - {\vec g}^h$. Algorithm~\ref{alg:relax} illustrates $\mu$ sweeps of Richardson iteration applied to (\ref{resid-eq}). 

\noindent
\begin{center}
\begin{minipage}{1.\linewidth}
\begin{algorithm}[H]
\caption{$\algG^h\left({\vec v}^h, {\vec g}^h,\mu,h\right)$; Relaxation}
\label{alg:relax}
\begin{algorithmic}[1]
  \Require ${\mat L}^h$
  \State $i \leftarrow 0$
  \While{$i<\mu$}
  \State ${\vec v}^h \leftarrow {\vec v}^h - \frac{1}{\|{\mat L}^h\|}\left({\mat L}^h{\vec v}^h - {\vec g}^h\right)$
  \COMMENT{Richardson iteration}
  \State $i \leftarrow i+1$  
  \EndWhile
  \State\Return ${\vec v}^h$
\end{algorithmic}
\end{algorithm}
\end{minipage}
\end{center}
\vspace{4mm}

To show that Richardson converges with the chosen scale, consider first the following simple relationship between the error and the residual in terms of the target equation in (\ref{resid-eq}):
\begin{equation}
{\mat L}^h {\vec e}^h = {\mat L}^h \left({\vec v}^h - ({\mat L}^h)^{-1}{\vec g}^h\right) = {\mat L}^h{\vec v}^h - {\vec g}^h = {\vec r}^h.
\label{res-err}
\end{equation}
If we use this relationship in (\ref{error}) and subtract $({\mat L}^h)^{-1}{\vec g}^h$ from both sides of the result, then with ${\mat I}^h \in \R^{n \times n}$ denoting the identity on grid $h$, we end up with the error propagation expression
\begin{equation}
{\vec e}^h \leftarrow {\mat G}^h{\vec e}^h, \mbox{   where  }{\mat G}^h = {\mat I}^h - \frac{1}{\|{\mat L}^h\|}{\mat L}^h.
\label{g}
\end{equation}
This ${\mat G}^h$ is called the error propagation matrix because it governs how the iteration transforms the error. To show that Richardson's iteration is convergent in the Euclidean norm $\|\cdot\|$, let $\langle \cdot, \cdot \rangle$ denote the Euclidean inner product and assume that ${\vec e}^h$ is any {\em nonzero} error component. Then
\begin{equation}
0 < \langle \frac{1}{\|{\mat L}^h\|}{\mat L}^h {\vec e}^h, {\vec e}^h \rangle \le 
\frac{\|{\mat L}^h\| \langle {\vec e}^h, {\vec e}^h \rangle}{ \|{\mat L}^h\|} = \langle {\vec e}^h, {\vec e}^h \rangle,
\label{ineq}
\end{equation}
which can be rewritten (from right to left) as 
\begin{equation}
\label{convergent}
0 \le \langle ({\mat I}^h - \frac{1}{\|{\mat L}^h\|}{\mat L}^h){\vec e}^h, {\vec e}^h \rangle < \langle {\mat I}^h {\vec e}^h, {\vec e}^h \rangle = \langle {\vec e}^h, {\vec e}^h \rangle.
\end{equation}
This bound means that $\|{\mat G}^h\| < 1$, which proves that Richardson's iteration is convergent in the Euclidean norm because it implies that the new error in (\ref{g}) must be smaller than the old error. Unfortunately, convergence does not necessarily mean fast convergence as we discuss next.

\subsection{Convergence behavior} Simple relaxation methods applied to discretizations of PDEs, including Richardson's iteration, typically stall well before they reach an acceptable approximation to the solution.  The cause of stalling comes from the residual's inability to see {\em algebraically smooth} error, by which we mean that the matrix applied to such an error yields a scaled residual $\frac{1}{\|{\mat L}^h\|}({\mat L}^h{\vec v}^h - {\vec g}^h)$ that is small compared to the {\em algebraic} error ${\vec e}^h \equiv {\vec v}^h - ({\mat L}^h)^{-1}{\vec g}^h$. This means that the scaled residual does little to improve ${\vec v}^h$. For many discretized elliptic equations, this algebraic sense of smoothness of the error usually corresponds to the geometric sense, where errors vary slowly across the grid. But we adopt the algebraic sense here instead because it applies to more general matrix equations that may not even have a geometric basis, and it sheds important light on how to treat these errors. 

It is unfortunate enough that correcting such algebraically smooth errors by way of the residual, as simple relaxation methods do, would accomplish very little. But worse yet, while these methods may work well for a couple of iterations when the initial error has oscillatory components, they shoot themselves in the foot because this fast elimination of oscillatory error exposes the remaining algebraically smooth error that works to stall all subsequent iterations. While this limitation is a common smoothing property of most conventional iterative methods applied to discretizations of PDEs, we show below that it is just what makes coarsening work.

\begin{figure*}[t!]
        \centering
        \subfloat[Initial guess]{
            \includegraphics[width=.30\linewidth]{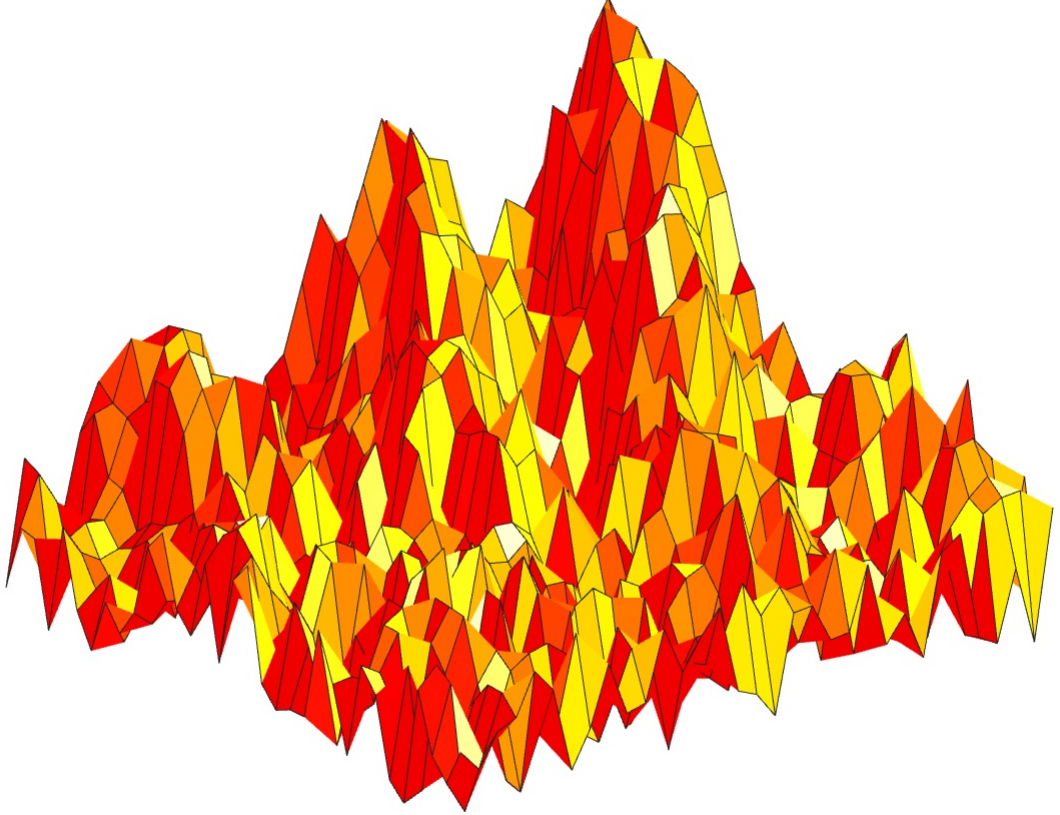}
        }
        \subfloat[Iteration 1]{
            \includegraphics[width=.30\linewidth]{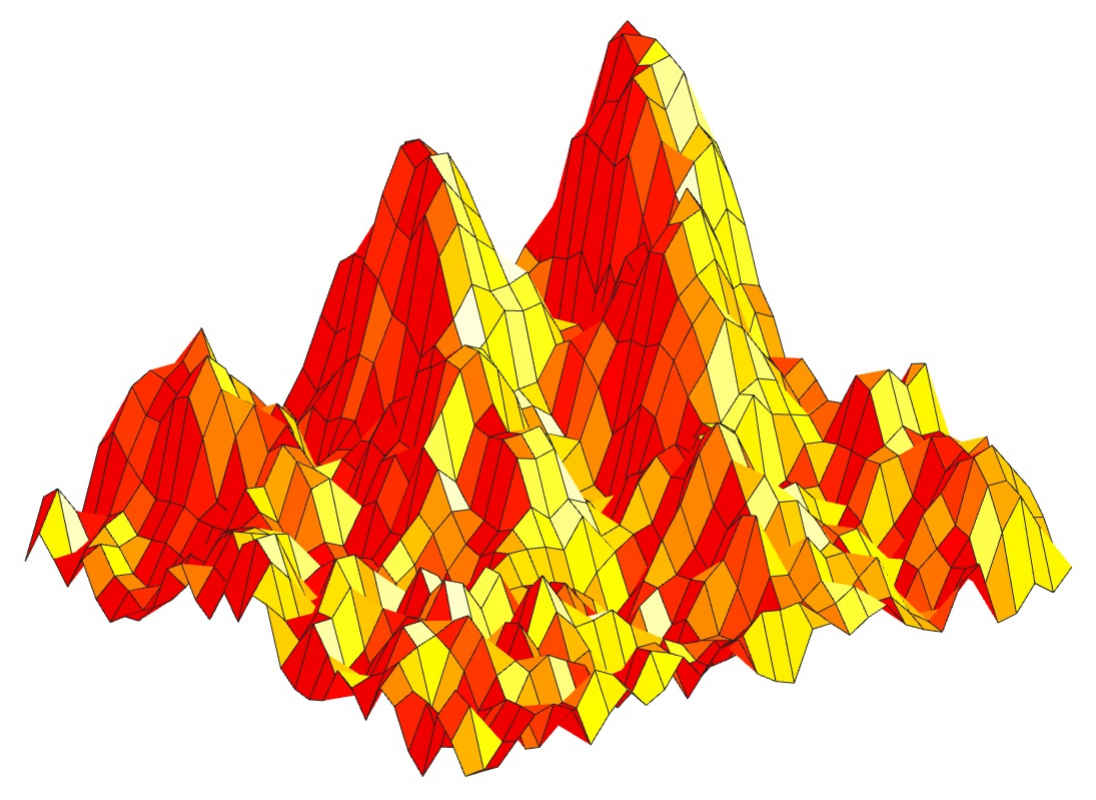}
        }
        \subfloat[Iteration 2]{
            \includegraphics[width=.30\linewidth]{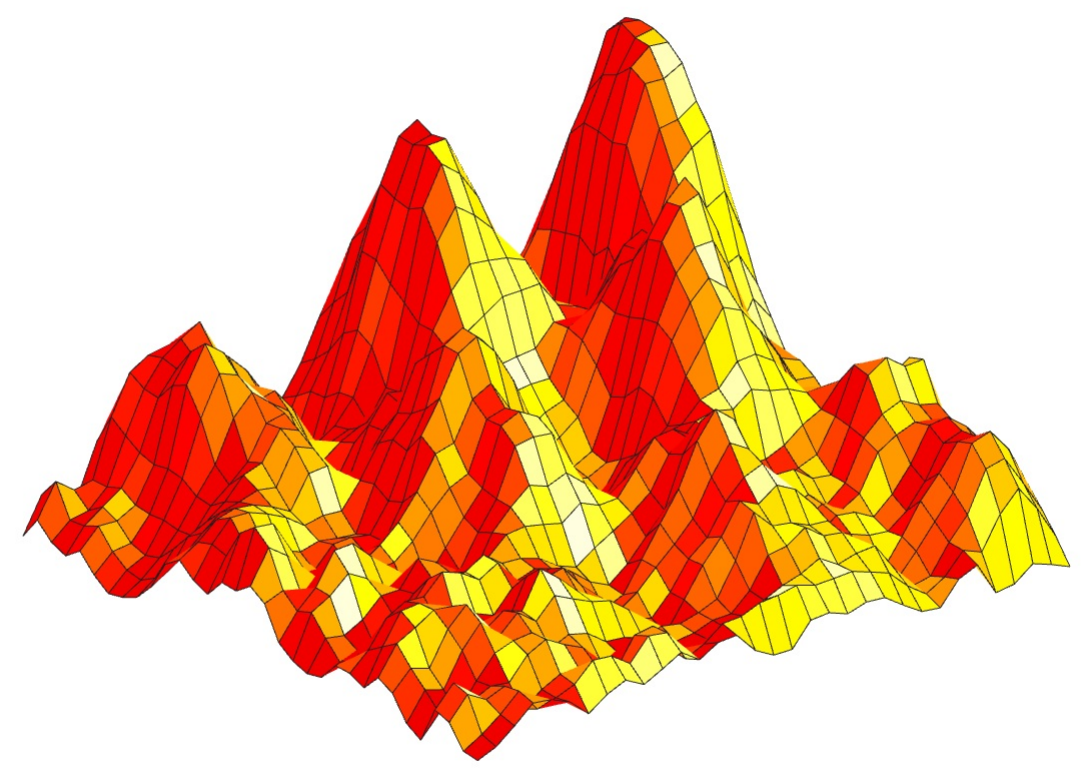}
        }\\
        \subfloat[Iteration 3]{
            \includegraphics[width=.30\linewidth]{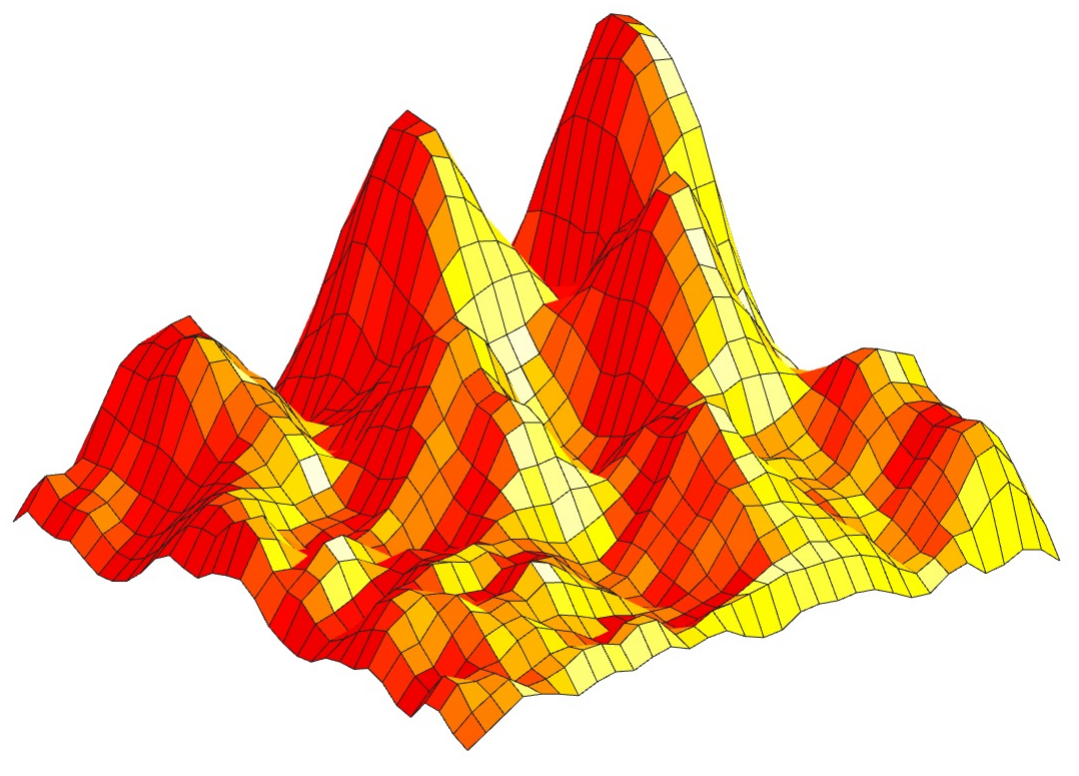}
        }
        \subfloat[Iteration 10]{
            \includegraphics[width=.30\linewidth]{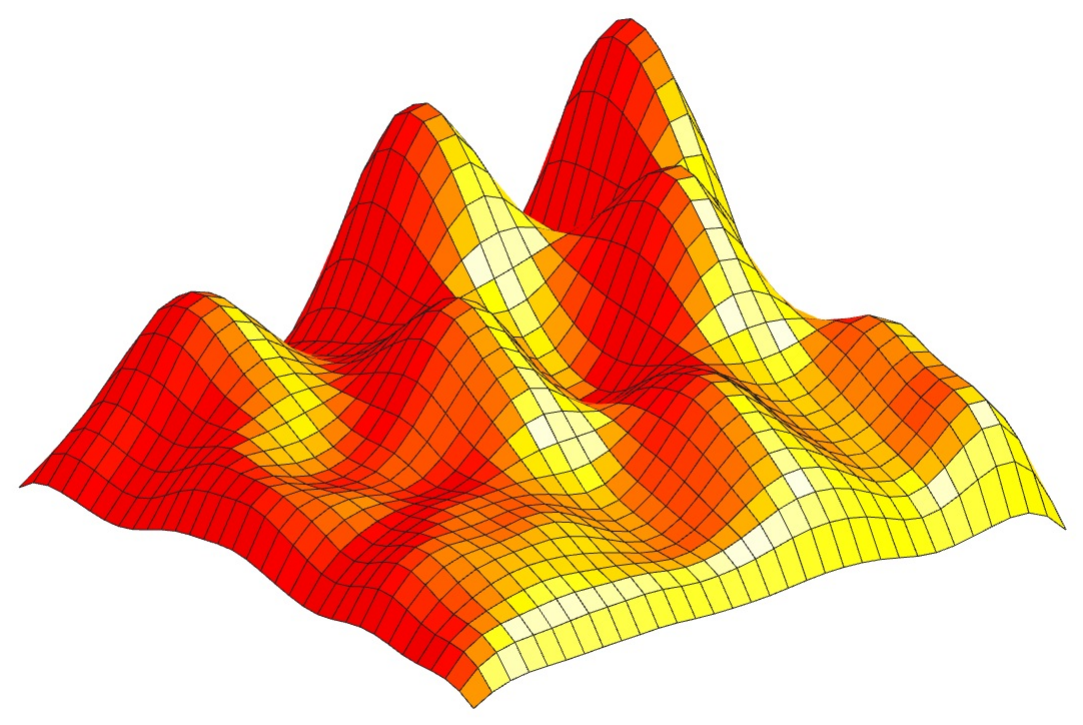}
        }
        \subfloat[Iteration 20]{
            \includegraphics[width=.30\linewidth]{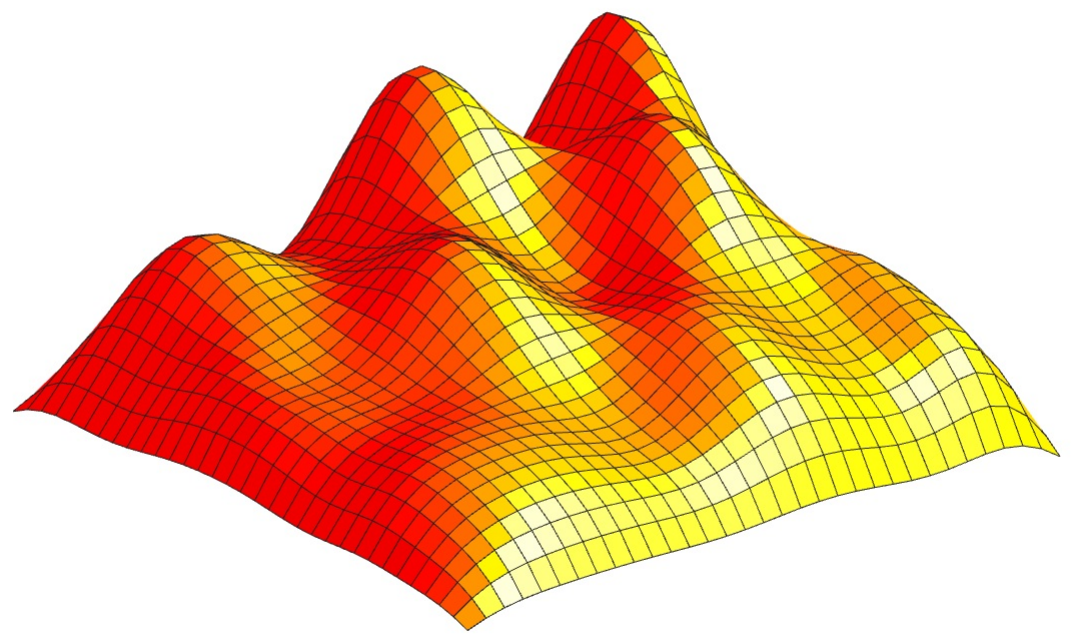}
        }
        \caption{Richardson relaxation applied to the Poisson problem on a unit square. What is shown is the magnitude of the error relative to a known reference solution. While oscillatory components of the error are quickly eliminated, smooth components remain even after many iterations.}
        \label{fig:fig}
    \end{figure*}

This smoothing behavior of Richardson iteration is illustrated in Figure~\ref{fig:fig}. The problem is (\ref{resid-eq}) with ${\vec g}^h = {\vec 0}^h$ and ${\mat A}^h$ as the so-called five-point discrete Poisson operator formed from continuous piecewise linear finite elements on a uniform grid of mesh spacing $h=\frac{1}{32}$ on the two-dimensional unit square with homogeneous boundary conditions.\footnote{In this model Poisson case, errors are smoothed in both the geometric and the algebraic senses.} The homogeneous assumptions mean that we are solving for ${\vec v}^h={\vec 0}^h$ in the equation ${\mat A}^h{\vec v}^h = {\vec 0}^h$, so the error in any approximation ${\vec v}^h$ is simply ${\vec v}^h$ itself. This tends to avoid any bothersome rounding errors, especially after many iterations, but its most important feature is knowing exactly what the error is. The initial guess is 
\[
\tfrac{1}{\alpha}\left(\sin(1.4x + 0.1)\sin(1.4y + 0.1)(1 + \sin(17x - 2)\sin(9y)) + 0.7\mathrm{rand}(0,1)\right),
\]
where $\alpha$ is chosen such that the maximum absolute value is $1$. This guess consists of a combination of two geometrically smooth errors and one geometrically oscillatory error, which offers a good visualization of how smoothing works. Note how the oscillatory component is fairly quickly reduced relative to the emerging smooth parts that remain. We study this slowness in some detail next.

\subsection{Algebraic smoothness measures}While (\ref{convergent}) shows that the error is reduced, it does not say anything about smoothness. For a given nonzero error ${\vec e}^h$, we therefore introduce the following two quantities that are particularly useful for measuring its algebraic smoothness:
$$\M_w({\vec e}^h) = \frac{\langle {\mat L}^h{\vec e}^h,{\vec e}^h \rangle }{||{\mat L}^h||\langle {\vec e}^h,{\vec e}^h \rangle}
\mbox{   and   } \M_s({\vec e}^h)= \frac{\langle {\mat L}^h{\vec e}^h,{\mat L}^h{\vec e}^h \rangle }{||{\mat L}^h||\langle {\mat L}^h{\vec e}^h,{\vec e}^h \rangle}.$$ 
The {\em weak} measure, $\M_w({\vec e}^h)$, is the relative Rayleigh quotient of ${\vec e}^h$, and the {\em strong} measure, $\M_s({\vec e}^h)$, is the relative Rayleigh quotient of $\left({\mat L}^h\right)^{\frac{1}{2}}{\vec e}^h$, both relative to $\|{\mat L}^h\|$. Because $\|{\mat L}^h\|$ is the maximum eigenvalue of ${\mat L}^h$, note that the ranges of both measures are in $(0,1]$. In what follows, we call errors with a large measure algebraically oscillatory or simply oscillatory, while errors with small measure are called algebraically smooth or simply smooth.

When ${\mat L}^h$ is symmetric, as we have assumed here, any error can be written as a linear combination of the eigenvectors, $\{{\vec w}_1,\dots,{\vec w}_n$\}, of ${\mat L}^h$: ${\vec e}^h = \sum_{i=1}^n \alpha_i {\vec w}_i$. Smooth vectors are rich in eigenvectors associated with the low part of the spectrum, while oscillatory vectors are rich in the high part. This is the basis for using the terms \emph{smooth} and \emph{oscillatory}. For example, if ${\vec e}^h$ is  an eigenvector of ${\mat L}^h$ associated with eigenvalue $\lambda$, then $\M_w({\vec e}^h)$ = $\M_s({\vec e}^h)$ = $\frac{\lambda}{\|{\mat L}^h\|}$. By virtue of our SPD assumption, $\|{\mat L}^h\| = \max_{1\le i \le n} \lambda_i,$ where the $\lambda_i$ are the eigenvalues of ${\mat L}^h.$ Thus, $\M_w({\vec e}^h)$ assesses how large $\lambda$ is relative to its largest possible value. For general linear combinations of eigenvectors, this measure assesses where the error lives in terms of the spectrum of ${\mat L}^h$. For many applications, this value can be very small at the lower end of the spectrum (typically $O(h^2)$ for discrete second-order elliptic operators). Conditions on the PDE and its discretization that produce matrices satisfying the strong approximation property are given, for example, in \cite{vassilevski2010}.

These weak and strong measures are important because they identify errors that cannot be adequately reduced by relaxation and therefore must be reduced by coarse-grid correction. To begin to see this, note from (\ref{g}) that relaxation reduces the nonzero current error, ${\vec e}^h$, in the Euclidean norm by the factor $\|{\mat G}^h{\vec e}^h\|/\|{\vec e}^h\|$. We have that
\begin{equation}
\frac{\|{\mat G}^h{\vec e}^h\|^2}{\|{\vec e}^h\|^2} = 1 - 2 \frac{\langle {\mat L}^h {\vec e}^h, {\vec e}^h \rangle }{\|{\mat L}^h\| \cdot \|{\vec e}^h\|^2} 
+ \frac{\langle {\mat L}^h {\vec e}^h, {\mat L}^h{\vec e}^h \rangle }{\|{\mat L}^h\|^2 \|{\vec e}^h\|^2} = 1 - (1 + \xi) \M_w({\vec e}^h),
\label{xi}
\end{equation}
where
\[
\xi = \frac{\langle {\mat L}^h {\vec e}^h, {\vec e}^h \rangle }{\|{\mat L}^h\| \cdot \|{\vec e}^h\|^2} 
- \frac{\langle {\mat L}^h {\vec e}^h, {\mat L}^h{\vec e}^h \rangle }{\|{\mat L}^h\|^2 \|{\vec e}^h\|^2} .
\]
Note that $\xi \in [0,1)$ because
\[
0 < \frac{\langle {\mat L}^h {\vec e}^h, {\mat L}^h{\vec e}^h \rangle }{\|{\mat L}^h\|^2 \|{\vec e}^h\|^2} 
\le \frac{\|{\mat L}^h\| \langle {\mat L}^h {\vec e}^h, {\vec e}^h \rangle }{\|{\mat L}^h\|^2 \cdot \|{\vec e}^h\|^2}
= \frac{\langle {\mat L}^h {\vec e}^h, {\vec e}^h \rangle }{\|{\mat L}^h\| \cdot \|{\vec e}^h\|^2} = \M_w({\vec e}^h) .
\]
Estimate (\ref{xi}) confirms that relaxation slows {\em in the Euclidean norm} if and only if the {\em weak measure} is small. This can be seen by rewriting (\ref{xi}):
\[
\M_w({\vec e}^h) = \frac{1-\frac{\|{\mat G}^h{\vec e}^h\|^2}{\|{\vec e}^h\|^2}}{1 + \xi} .
\]

Let $\langle \cdot,\cdot\rangle_{{\mat L}^h} = \langle {\mat L}^h\cdot,\cdot\rangle $ and $\|\cdot\|_{{\mat L}^h} \equiv \sqrt{\langle \cdot,\cdot\rangle_{{\mat L}^h}}$ denote the {\em energy} inner product and its induced {\em energy} norm, respectively. Then an argument similar to the above shows that
\[
\frac{\|{\mat G}^h{\vec e}^h\|_{{\mat L}^h}^2}{\|{\vec e}^h\|_{{\mat L}^h}^2} = 1 - (1 + \chi) \M_s({\vec e}^h),
\]
for some $\chi \in [0,1)$. This expression shows that relaxation converges in the energy norm because the right-hand side is less than $1$. But it also shows that relaxation slows in energy if and only if the {\em strong measure} is small. 

This correspondence between weak vs. strong measures and Euclidean vs. energy norms carries over to the analysis of coarse-grid correction. As shown in \cite{Vassilevski2008}, if the coarse grid adequately approximates errors for which $\M_s({\vec e}^h)$ is small, then the so-called multigrid {\em V-cycles} converge well in the energy norm. For the so-called {\em two-grid} or {\em W-cycles} (described in what follows) in the Euclidean norm, it is enough to have the coarse grid adequately approximate errors for which $\M_w({\vec e}^h)$ is small. This latter requirement is weaker in part because small $\M_s({\vec e}^h)$ implies small $\M_w({\vec e}^h)$, which follows because
\[
\langle {\mat L}^h{\vec e}^h,{\vec e}^h \rangle^2 \le 
\| {\mat L}^h{\vec e}^h \|^2 \|{\vec e}^h\|^2 =
\langle {\mat L}^h{\vec e}^h,{\mat L}^h{\vec e}^h \rangle\langle {\vec e}^h,{\vec e}^h \rangle .
\]
The converse, however, is not true, which means that aiming for $\M_s({\vec e}^h)$ to be small is the more restrictive objective, hence the term {\em strong}. More importantly, what we meant by adequate approximation for $\M_s({\vec e}^h)$ involves the energy norm, which is a stronger requirement than that for the Euclidean norm used with the weak measure.

\subsection{Other relaxation methods} We have relied on the simplicity of Richardson, Chebyshev, and other Krylov methods, where error propagation matrices (as defined in (\ref{g})) are polynomials in ${\mat L}^h$. For other choices, the analysis and attendant algebra can be much more involved, although most of the basic principles are essentially the same. It is also important to point out that dividing the residual by $\|{\mat L}^h\|$ ensures that the iteration converges, but it is not necessarily the best choice for relaxation. See, for example, \cite{briggs00} for somewhat better scale factors. 

Krylov methods by themselves may not be effective at smoothing the error in a geometric sense even in a fairly simple setting. One limitation is that is uses a residual that does not account for variations in scale across the grid. For instance, beyond the simple uniform grid setting treated here, suppose for instance that the matrix has a block, corresponding to a local subdomain of the PDE, of very large entries relative to the rest of the matrix entries. This disparity could be due to a relatively fine resolution there or relatively large associated coefficients of the PDE. Then Richardson would quickly reduce the oscillatory error in that block but cannot effectively reduce the oscillatory error elsewhere because the associated residuals would be too small. A simple solution to this difficulty is to use Jacobi's method given by
\[
{\vec v}^h \leftarrow {\vec v}^h - s \left({\mat D}^h\right)^{-1}\left({\mat L}^h{\vec v}^h - {\vec g}^h\right),
\]
where ${\mat D}^h$ is the diagonal part of ${\mat L}^h$ and $s>0$ is a suitably chosen step size.

In any event with the understanding of the smoothness property of Richardson in hand, we now turn to coarsening, the other basic multigrid component.

\section{Coarsening}
\label{sec:coarseCorrection}

The coarse-grid correction phase of multigrid enters as a way to exploit the smooth errors that relaxation tends to produce. The basic idea is that algebraically smooth error varies so slowly from one neighboring grid point to the next that it can be adequately approximated by fewer grid points. The error components thus computed on these coarse levels are then interpolated back to the fine grid to correct the approximation there. We see how this could be done next by developing a structure of coarsening based on an energy principle.

\subsection{Energy functional} Consider the grid $h$ {\em energy functional} given by
\[
\F^h({\vec v}^h) \equiv \langle {\mat L}^h{\vec v}^h,{\vec v}^h\rangle - 2 \langle {\vec v}^h,{\vec g}^h\rangle .
\]
To show that the minimum of $\F^h({\vec v}^h)$ is the solution of (\ref{resid-eq}), let $C \equiv \langle  {\vec g}^h,({\mat L}^h)^{-1}{\vec g}^h\rangle$ and note that 
\begin{align*}
\F^h({\vec v}^h) &= \langle {\mat L}^h{\vec v}^h,{\vec v}^h\rangle - 2 \langle {\vec v}^h,{\vec g}^h\rangle + \langle  {\vec g}^h,({\mat L}^h)^{-1}{\vec g}^h\rangle -  \langle  {\vec g}^h,({\mat L}^h)^{-1}{\vec g}^h\rangle \\
&= \langle {\mat L}^h{\vec v}^h,{\vec v}^h\rangle - 2 \langle {\mat L}^h{\vec v}^h,({\mat L}^h)^{-1}{\vec g}^h\rangle + \langle {\mat L}^h ({\mat L}^h)^{-1}{\vec g}^h,({\mat L}^h)^{-1}{\vec g}^h\rangle -  C \\
&= \langle{\mat L}^h\left({\vec v}^h - ({\mat L}^h)^{-1}{\vec g}^h\right), {\vec v}^h - ({\mat L}^h)^{-1}{\vec g}^h\rangle -  C \\
&= \|{\vec e}^h\|_{{\mat L}^h}^2 -  C .
\end{align*}
Thus, $\F^h({\vec v}^h)$ differs from the squared energy norm of the error by a constant ($C$ does not depend on ${\vec v}^h$), which means that ${\vec v}^h$ solves (\ref{resid-eq}) if and only if it is the minimizer of $\F^h({\vec v}^h)$. (That minimum is unique because the above equation shows that the minimum value of $-C$ is attained only when ${\vec e}^h \equiv {\vec v}^h -({\mat L}^h)^{-1}{\vec g}^h = {\vec 0}^h$.) 

Assume, as before, that grid $2h$ and an interpolation operator ${\mat P}_{2h}^h$ from grid $2h$ to grid $h$  are available. Assume also that ${\mat P}_{2h}^h$ is of full rank to avoid singularity in the coarse-grid matrix constructed below. In any event, we can now exploit this minimization principle further by noting that the best coarse-grid correction to a fixed approximation, ${\vec v}^h$, in the sense of minimizing $\F^h({\vec v}^h - {\mat P}_{2h}^h {\vec v}^{2h})$ over all possible ${\vec v}^{2h}$, is expressed by
\begin{equation}
{\vec v}^h \leftarrow {\vec v}^h - {\mat P}_{2h}^h\left(\left({\mat P}_{2h}^h\right)^t {\mat L}^h{\mat P}_{2h}^h\right)^{-1}\left({\mat P}_{2h}^h\right)^t({\mat L}^h{\vec v}^h - {\vec g}^h).
\label{correction}
\end{equation}
(We use superscript $t$ to denote matrix transpose.) To verify this form of the correction, note that
\[
\langle {\mat L}^h {\mat P}_{2h}^h{\vec v}^{2h}, {\vec v}^h\rangle = \langle {\vec v}^{2h}, \left({\mat P}_{2h}^h\right)^t{\mat L}^h {\vec v}^h\rangle ,
\]
which shows that
\begin{align*}
&\F^h({\vec v}^h - {\mat P}_{2h}^h {\vec v}^{2h}) \\
&= \langle {\mat L}^h\left({\vec v}^h  - {\mat P}_{2h}^h{\vec v}^{2h}\right),{\vec v}^h  - {\mat P}_{2h}^h{\vec v}^{2h}\rangle - 2 \langle {\vec v}^h  - {\mat P}_{2h}^h{\vec v}^{2h}),{\vec g}^h\rangle \\
&= \langle {\mat L}^h{\vec v}^h,{\vec v}^h\rangle - 2 \langle {\mat L}^h{\mat P}_{2h}^h{\vec v}^{2h}, {\vec v}^h\rangle + \langle {\mat L}^h{\mat P}_{2h}^h{\vec v}^{2h},{\mat P}_{2h}^h{\vec v}^{2h}\rangle - 2 \langle {\vec v}^h,{\vec g}^h\rangle + 2 \langle {\mat P}_{2h}^h{\vec v}^{2h},{\vec g}^h\rangle \\
&= \langle {\mat L}^h{\vec v}^h,{\vec v}^h\rangle - 2 \langle {\vec v}^h,{\vec g}^h\rangle + \langle {\mat L}^h{\mat P}_{2h}^h{\vec v}^{2h},{\mat P}_{2h}^h{\vec v}^{2h}\rangle - 2 \langle {\mat L}^h{\mat P}_{2h}^h{\vec v}^{2h}, {\vec v}^h\rangle + 2 \langle {\mat P}_{2h}^h{\vec v}^{2h},{\vec g}^h\rangle \\
&= \langle {\mat L}^h{\vec v}^h,{\vec v}^h\rangle - 2 \langle {\vec v}^h,{\vec g}^h\rangle + \langle {\mat L}^h{\mat P}_{2h}^h{\vec v}^{2h},{\mat P}_{2h}^h{\vec v}^{2h}\rangle - 2 \langle {\vec v}^{2h}, \left({\mat P}_{2h}^h\right)^t({\mat L}^h{\vec v}^h - {\vec g}^h)\rangle \\
&= \F^h({\vec v}^h) + \F^{2h}({\vec v}^{2h}),
\end{align*}
where $\F^{2h}$ is the following grid $2h$ version of $\F^h$:
\[
\F^{2h}({\vec v}^{2h}) = \langle \left({\mat P}_{2h}^h\right)^t{\mat L}^h{\mat P}_{2h}^h{\vec v}^{2h},{\vec v}^{2h}\rangle - 2 \langle {\vec v}^{2h}, \left({\mat P}_{2h}^h\right)^t({\mat L}^h{\vec v}^h - {\vec g}^h)\rangle.
\]

This result implies that ${\vec v}^{2h}$ minimizes $\F^h({\vec v}^h - {\mat P}_{2h}^h {\vec v}^{2h})$ if and only if it minimizes $\F^{2h}({\vec v}^{2h})$. This grid $2h$ energy functional corresponds to the following grid $2h$ version of the matrix equation in (\ref{resid-eq}):
\begin{equation}
{\mat L}^{2h}{\vec v}^{2h} = {\vec g}^{2h},
\label{equation2h}
\end{equation}
where
\[
{\mat L}^{2h} = \left({\mat P}_{2h}^h\right)^t{\mat L}^h{\mat P}_{2h}^h \textrm{  and  } {\vec g}^{2h} = \left({\mat P}_{2h}^h\right)^t({\mat L}^h{\vec v}^h - {\vec g}^h) .
\]
Thus, the grid $2h$ correction is given by (\ref{correction}) as asserted.

\subsection{Variational conditions}This is the form of the coarse-grid correction that we use here, and it gives rise to the following so-called {\em variational conditions}:
\[
{\mat R}_h^{2h} \equiv \alpha^h \left({\mat P}_{2h}^h\right)^t  \mbox{   and   } {\mat L}^{2h} \equiv {\mat R}_h^{2h} {\mat L}^h{\mat P}_{2h}^h,
\]
where ${\mat R}_h^{2h}$ is the transfer operator from the fine to the coarse grid and ${\mat L}^{2h}$ is the coarse-grid matrix. The second definition is called the {\em Galerkin condition}, which specifies the coarse-grid operator. The first allows for any desired scaling $\alpha^h > 0$ of restriction; it does not contradict the variational principle because the correction in (\ref{correction}) is independent of $\alpha^h$ as the following indicates:
\[
\left(\left(\alpha^h{\mat P}_{2h}^h\right)^t {\mat L}^h{\mat P}_{2h}^h\right)^{-1}\left(\alpha^h{\mat P}_{2h}^h\right)^t =
\frac{1}{\alpha^h} \alpha^h \left(\left({\mat P}_{2h}^h\right)^t {\mat L}^h{\mat P}_{2h}^h\right)^{-1}\left({\mat P}_{2h}^h\right)^t.
\]
This invariance allows for a little flexibility in choosing these transfer operators so that they approximately preserve constants, for example. 

Because restriction and the coarse-grid problem follow directly from the choice of interpolation, the only basic choice to be made in designing a {\em variational} multigrid method besides relaxation is ${\mat P}_{2h}^h$. This choice amounts to determining a coarse set of grid points (or variables) and a method for interpolating functions defined on those points to functions defined on the fine-grid points. This energy-minimization formulation is convenient in the sense that, once this is done, the coarse-grid matrix and the restriction operator are then determined automatically from the variational conditions. This simplification is a major advantage that the energy minimization principle provides because choosing restriction separately is not always straightforward: choosing ${\mat R}_h^{2h}$ and ${\mat P}_{2h}^h$ to work together requires careful coordination and it means that an error measure may not be immediately at hand.

Having the structure for coarse-grid correction in place, we now introduce multigrid cycling. We start with the case of two grids as the foundation for more levels.

\section{Two-grid cycles}
\label{sec:2grid}

To apply the formal definition of the two-grid method introduced in the Introduction to the energy functional $\F^h({\vec v}^h)$, we first need to compute its derivative $\F^h({\vec v}^h)$ in the direction of ${\vec w}^h$. Remembering that the residual is ${\vec r}^h = {\mat L}^h{\vec v}^h - {\vec g}^h$, we have that
\begin{align*}
\F^h({\vec v}^h &+ s{\vec w}^h)) - \F^h({\vec v}^h) \\
&=\langle {\mat L}^h({\vec v}^h+s{\vec w}^h),{\vec v}^h+s{\vec w}^h\rangle - 2 \langle {\vec v}^h+s{\vec w}^h,{\vec g}^h\rangle - \langle {\mat L}^h{\vec v}^h,{\vec v}^h\rangle - 2 \langle {\vec v}^h,{\vec g}^h\rangle \\
&= \F^h({\vec v}^h) + 2s \left(\langle {\mat L}^h{\vec w}^h,{\vec v}^h\rangle - \langle {\vec w}^h,{\vec g}^h\rangle\right)
+ s^2 \langle {\mat L}^h{\vec w}^h,{\vec w}^h\rangle - \F^h({\vec v}^h)\\
&= 2s \left(\langle {\vec w}^h,{\mat L}^h{\vec v}^h\rangle - \langle {\vec w}^h,{\vec g}^h\rangle\right)
+ s^2 \langle {\mat L}^h{\vec w}^h,{\vec w}^h\rangle \\
&= 2s \langle {\vec w}^h,{\vec r}^h\rangle
+ s^2 \langle {\mat L}^h{\vec w}^h,{\vec w}^h\rangle.
\end{align*}
This shows that the derivative that we seek is given by
\[
{\F^h}'({\vec v}^h){\vec w}^h = \lim_{s \rightarrow 0} \frac{1}{s} \left(\F^h({\vec v}^h + s{\vec w}^h) - {\F^h}'({\vec v}^h)\right)
= \langle {\vec w}^h,2{\vec r}^h\rangle,
\]
which in turn shows that the gradient of $\F^h({\vec v}^h)$ is $2{\vec r}^h$. Absorbing the $2$ into $t^h$, then the formal two-grid method becomes
\vspace{1.0em}
\begin{itemize}
\item Set ${\vec v}^h \leftarrow {\vec v}^h - s^h{\vec r}^h$, where $s^h$ minimizes $\F^h({\vec v}^h - t^h{\vec r}^h)$ over all $t^h \in \Re$. 
\item Set ${\vec v}^h \leftarrow {\vec v}^h - {\mat P}_{2h}^h {\vec v}^{2h}$, where ${\vec v}^{2h}$ minimizes $\F^h({\vec v}^h - {\mat P}_{2h}^h {\vec w}^{2h})$ over all ${\vec w}^{2h}$ on grid $2h$.
\end{itemize}
\vspace{01.0em}

To make this formal definition concrete, consider the following specific two-grid method, as illustrated in Algorithm~\ref{alg:2-grid}:
\vspace{1.0em}
\begin{itemize}
\item Relax on ${\vec v}^h$: apply ${\vec v}^h \leftarrow {\vec v}^h - \frac{1}{\|{\mat L}^h\|}\left({\mat L}^h{\vec v}^h - {\vec g}^h\right)$ $\mu$ times.
\item Apply coarse-grid correction: ${\vec v}^h \leftarrow {\vec v}^h - {\mat P}_{2h}^h ({\mat L}^{2h} )^{-1} {\mat R}_h^{2h} \left({\mat L}^h {\vec v}^h - {\vec g}^h\right).$
\item Relax on ${\vec v}^h$: apply ${\vec v}^h \leftarrow {\vec v}^h - \frac{1}{\|{\mat L}^h\|}\left({\mat L}^h{\vec v}^h - {\vec g}^h\right)$ $\nu$ times.
\end{itemize}
\vspace{1.0em}
Note that we have not use the {\em optimal} scalar\footnote{The optimal choice is $s^h = \langle {\vec r}^h, {\vec r}^h\rangle/\langle {\mat A}^h{\vec r}^h, {\vec r}^h\rangle$, which can be derived by finding the zero of the derivative  w.r.t. $t^h$ of $\F^h({\vec v}^h - t^h{\vec r}^h)$.} in the first step above but focussed instead on Richardson with the simpler choice $s^h = \frac{1}{\|{\mat L}^h\|}$. While not necessarily optimal, this choice does effectively smooth the error. On the other hand, we have included the optimal coarse-grid correction in the second step above because it merely states that ${\vec v}^{2h} = ({\mat L}^{2h} )^{-1} {\mat R}_h^{2h} \left({\mat L}^h {\vec v}^h - {\vec g}^h\right)$. Note also that we have incorporated $\mu$ pre-relaxation sweeps and $\nu$ post-relaxation sweeps, which is typical in practice: it is usually more effective to spread $\mu + \nu$ sweeps out before and after coarsening rather than bundling them up before or after because too many sweeps at once begin to struggle with smooth error that starts to dominate. Note, however, that in between the beginning and end of a sequence of two-grid cycles, there would be no difference between these choices in the two-grid case.

\noindent
\begin{center}
\begin{minipage}{1.\linewidth}
\begin{algorithm}[H]
\caption{$\algTG^h\left({\vec v}^h, {\vec g}^h,\mu, \nu, h\right)$; Two-grid cycle}
\label{alg:2-grid}
\begin{algorithmic}[1]
      \Require ${\mat L}^h$, ${\mat L}^{2h}$, ${\mat R}_h^{2h}$, ${\mat P}_{2h}^h$
      \State ${\vec v}^h \leftarrow \algG^h\left({\vec v}^h, {\vec g}^h,\mu,h\right)$
      \COMMENT{Relax $\mu$ times.}
      \If{$h<H$} 
      \COMMENT{Keep coarsening.}
      \State ${\vec g}^{2h} = {\mat R}_h^{2h} \left({\mat L}^h {\vec v}^h - {\vec g}^h\right)$
      \COMMENT{Transfer grid $h$ residual to grid {2h}.}
      \State ${\vec v}^{2h} \leftarrow \left({\mat L}^{2h}\right)^{-1}{\vec g}^{2h}$
      \COMMENT{Solve the grid $2h$ equation exactly.}
      \State ${\vec v}^h \leftarrow {\vec v}^h - {\mat P}_{2h}^h{\vec v}^{2h}$
      \COMMENT{Correct fine-grid iterate.}
      \State ${\vec v}^h \leftarrow \algG^h\left({\vec v}^h, {\vec g}^h,\nu,h\right)$
      \COMMENT{Relax $\nu$ times.}
      \EndIf
     \State\Return ${\vec v}^h$
\end{algorithmic}
\end{algorithm}
\end{minipage}
\end{center}
\vspace{4mm}

\begin{rem} {\bf A word about residuals.}  Our next task is to extend this two-grid method to multiple levels. To do so, the key is to recognize that the grid $2h$ equation (\ref{equation2h}) is of the same form as the grid $h$ equation (\ref{resid-eq}). This suggests that all we have to do is apply the two-grid cycle to grid $2h$ to get started with three levels that includes grid $4h$. This approach provides a natural way to recursively define multilevel cycles, but it also introduces potential confusion. Grid $h$ relaxation applies directly to (\ref{resid-eq}), but (\ref{equation2h}) is indirectly applied to finding a correction that solves an averaged residual equation for (\ref{resid-eq}). This recursion means that grid $2h$ ends up passing an averaged residual for its averaged residual equation to grid $4h$. We are talking here about averaged residuals of averaged residual equations and so on down the grid hierarchy. None but the finest grid involves the original source term ${\vec g}^h$ directly. This is why we have  introduced (\ref{resid-eq}) here to allow the source term to change within the multigrid cycle, while the notation in (\ref{equation}) is fixed. It is perhaps advisable for the reader to begin by thinking carefully about how multigrid works with just two grids, hopefully to better understand its recursive extension to many levels in what follows. \end{rem}

We appeal  next to recursion to extend the two-grid cycle to a full hierarchy of grids.

\section{Multilevel cycles}
\label{sec:solvers}

Two-grid versions usually only achieve true multigrid efficiency when grid $h$ is very coarse and therefore relatively easy to solve by relaxation alone. However, in this case, grid $h$ is also probably fairly easy to solve by relaxation alone, so multigrid would hardly be essential. The point here is that a fairly fine grid with several levels in the hierarchy allows multigrid to realize its full potential. The key to the extension to a potentially fast multilevel version is to notice that coarse-grid correction involves the solution of ${\mat L}^{2h}{\vec v}^{2h} = {\mat R}_h^{2h} \left({\mat L}^h{\vec v}^h - {\vec g}^h\right)$ (followed by correction to the fine grid via ${\vec v}^h \leftarrow {\vec v}^h - {\mat P}_{2h}^h {\vec v}^{2h}$). Given a hierarchy of coarse grids from $h$ to $2h$ to $4h$ and on down to a very coarse grid $H$, the idea is to simply replace this exact solution step by a scheme that first improves the grid $2h$ initial approximation ${\vec v}^{2h} = {\vec 0}^{2h}$ by $\mu$ pre-relaxation sweeps, a subsequent correction from grid $4h$, and $\nu$ post-relaxation sweeps on grid $2h$. 

Continuing recursively, this process proceeds in this way with $\mu$ sweeps on each level down to the coarsest grid $H$ and back up to the finest level with $\nu$ sweeps on each level. Note that the grid $H$ equation involves $\mu + \nu$ relaxation sweeps, although it is often the case in practice that relaxation there is replaced by a more aggressive method, such as Gaussian elimination. This would be especially important when grid $H$ cannot be easily chosen coarse enough for relaxation to converge quickly, but still coarse enough that such methods are feasible and relatively inexpensive. However, for very coarse grids, relaxation alone should suffice.

This multilevel process is called a $V(\mu,\nu)$-cycle because it starts from the fine grid and passes residuals down through the coarser grids to the coarsest, and then proceeds back to the finest by passing up corrections. This $V(\mu,\nu)$ cycle applied to (\ref{resid-eq}) is represented by the following expression: 
\begin{equation}
{\vec v}^h \leftarrow \algMV^h\left({\vec v}^h, {\vec g}^h,\mu,\nu,h\right),
\end{equation}
where the ${\vec v}^h$ on the right is a current grid $h$ approximation and ${\vec g}^h$ is a given source term. $\algMV^h\left({\vec v}^h, {\vec g}^h,\mu,\nu,h\right)$ is defined recursively in Algorithm~\ref{alg:V-cycle}.

\noindent
\begin{center}
\begin{minipage}{1.\linewidth}
\begin{algorithm}[H]
\caption{$\algMV^h\left({\vec v}^h, {\vec g}^h,\mu,\nu,h\right)$; V-cycle}
\label{alg:V-cycle}
\begin{algorithmic}[1]
      \Require ${\mat L}^h$, ${\mat R}_h^{2h}$, ${\mat P}_{2h}^h$
      \State ${\vec v}^h \leftarrow \algG^h\left({\vec v}^h, {\vec g}^h,\mu,h\right)$
      \COMMENT{Relax $\mu$ times.}
      \If{$h<H$} 
      \COMMENT{Keep coarsening.}
      \State ${\vec g}^{2h} = {\mat R}_h^{2h} \left({\mat L}^h {\vec v}^h - {\vec g}^h\right)$
      \COMMENT{Transfer grid $h$ residual to grid {2h}.}
      \State ${\vec v}^{2h} \leftarrow \algMV^{2h}\left({\vec 0}^{2h}, {\vec g}^{2h}, \mu,\nu,2h\right)$
      \COMMENT{Apply $\algMV^{2h}$ with zero initial guess.}
      \State ${\vec v}^h \leftarrow {\vec v}^h - {\mat P}_{2h}^h{\vec v}^{2h}$
      \COMMENT{Correct fine-grid iterate.}
      \EndIf
      \State ${\vec v}^h \leftarrow \algG^h\left({\vec v}^h, {\vec g}^h,\nu,h\right)$
      \COMMENT{Relax $\nu$ times.}
     \State\Return ${\vec v}^h$
\end{algorithmic}
\end{algorithm}
\end{minipage}
\end{center}
\vspace{4mm}

In some cases in practice, the approximation from the coarse grid is not accurate enough to achieve uniform V-cycle convergence. For example, piecewise constant interpolation used in coarsening approximates smooth fine-grid components too poorly to allow V-cycles to work effectively. For such cases, more work in the form of a stronger cycling scheme would be required, which brings us to the $W(\mu,\nu)$-cycle that coarsens twice between relaxation steps. It is represented by the expression
\begin{equation}
{\vec v}^h \leftarrow \algMW^h\left({\vec v}^h, {\vec g}^h,\mu,\nu,h\right)
\end{equation}
and defined recursively in Algorithm~\ref{alg:W-cycle}. 

\noindent
\begin{center}
\begin{minipage}{1.\linewidth}
\begin{algorithm}[H]
\caption{$\algMW^h\left({\vec v}^h, {\vec g}^h,\mu,\nu,h\right)$; $W$-cycle}
\label{alg:W-cycle}
\begin{algorithmic}[1]
      \Require ${\mat L}^h$, ${\mat R}_h^{2h}$, ${\mat P}_{2h}^h$
      \State ${\vec v}^h \leftarrow \algG^h\left({\vec v}^h, {\vec g}^h,\mu,h\right)$
      \COMMENT{Relax $\mu$ times.}
      \If{$h<H$} 
      \COMMENT{Keep coarsening.}
      \State ${\vec g}^{2h} = {\mat R}_h^{2h} \left({\mat L}^h {\vec v}^h - {\vec g}^h\right)$
      \COMMENT{Transfer grid $h$ residual to grid {2h}.}
      \State ${\vec v}^{2h} \leftarrow \algMW^{2h}\left({\vec 0}^{2h}, {\vec g}^{2h}, \mu,\nu,2h\right)$
      \COMMENT{Apply $\algMW^{2h}$ with zero initial guess.}
      \State ${\vec v}^{2h} \leftarrow \algMW^{2h}\left({\vec v}^{2h}, {\vec g}^{2h}, \mu,\nu,2h\right)$
      \COMMENT{Apply $\algMW^{2h}$ to grid $2h$ iterate.}
      \State ${\vec v}^h \leftarrow {\vec v}^h - {\mat P}_{2h}^h{\vec v}^{2h}$
      \COMMENT{Correct fine-grid iterate.}
      \EndIf
      \State ${\vec v}^h \leftarrow \algG^h\left({\vec v}^h, {\vec g}^h,\nu,h\right)$
      \COMMENT{Relax $\nu$ times.}
     \State\Return ${\vec v}^h$
\end{algorithmic}
\end{algorithm}
\end{minipage}
\end{center}
\vspace{4mm}

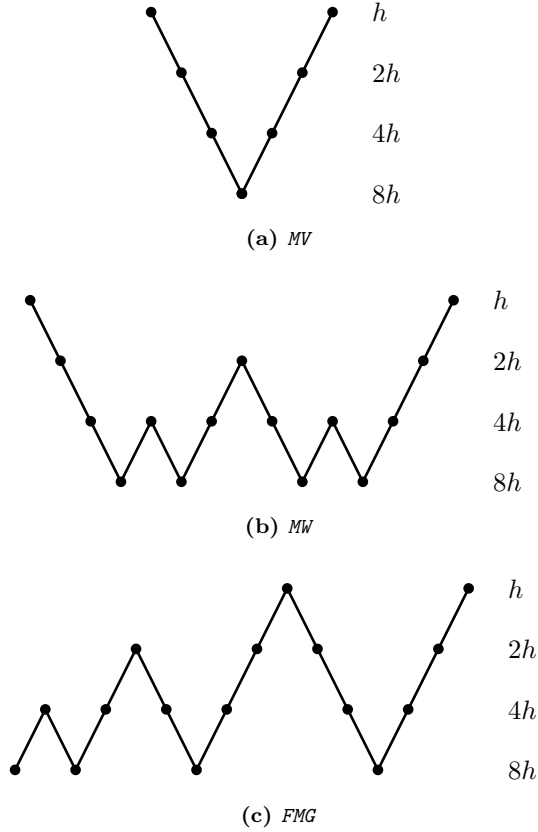
\begin{figure}[htbp]
\centering

\begin{subfigure}[t]{\textwidth}
\centering
\begin{tikzpicture}[scale=0.8]
  \draw[line width=1pt] (0,3) -- (1.5,0) -- (3,3);

  \foreach \y in {0,1,2,3} {
    \path (0,3) -- (1.5,0) coordinate[pos=\y/3] (L\y);
    \fill (L\y) circle (2.5pt);
  }

  \foreach \y in {0,1,2,3} {
    \path (3,3) -- (1.5,0) coordinate[pos=\y/3] (R\y);
    \fill (R\y) circle (2.5pt);
  }

  \node[right] at (3.5,3) {$h$};
  \node[right] at (3.5,2) {$2h$};
  \node[right] at (3.5,1) {$4h$};
  \node[right] at (3.5,0) {$8h$};
\end{tikzpicture}
\caption{$\algMV$}
\end{subfigure}


\begin{subfigure}[t]{\textwidth}
\centering
\begin{tikzpicture}[scale=0.8]
  \draw[line width=1pt] (0,3) -- (1.5,0) -- (2,1) -- (2.5,0) -- (3.5,2) -- (4.5,0) -- (5,1) -- (5.5,0) -- (7,3);

  \foreach \y in {0,1,2,3} {
    \path (0,3) -- (1.5,0) coordinate[pos=\y/3] (L\y);
    \fill (L\y) circle (2.5pt);
  }

  \foreach \y in {0,1,2} {
    \path (3.5,2) -- (2.5,0) coordinate[pos=\y/2] (R\y);
    \fill (R\y) circle (2.5pt);
  }
  
  \foreach \y in {1,2} {
    \path (3.5,2) -- (4.5,0) coordinate[pos=\y/2] (L\y);
    \fill (L\y) circle (2.5pt);
  }

  \foreach \y in {0,1,2,3} {
    \path (7,3) -- (5.5,0) coordinate[pos=\y/3] (R\y);
    \fill (R\y) circle (2.5pt);
  }

  \fill (2,1) circle (2.5pt);
  \fill (5,1) circle (2.5pt); 

  \node[right] at (7.5,3) {$h$};
  \node[right] at (7.5,2) {$2h$};
  \node[right] at (7.5,1) {$4h$};
  \node[right] at (7.5,0) {$8h$};
\end{tikzpicture}
\caption{$\algMW$}
\end{subfigure}

\begin{subfigure}[t]{\textwidth}
\centering
\begin{tikzpicture}[scale=0.8]
  \draw[line width=1pt] (0,0) -- (0.5,1) -- (1,0) -- (1.5,1) -- (2,2) -- (2.5,1) -- (3,0) -- (3.5,1) -- (4,2) -- (4.5,3) -- (5,2) -- (5.5,1) -- (6,0) -- (6.5,1) -- (7,2) -- (7.5,3);

  \fill (0,0) circle (2.5pt);
  \fill (0.5,1) circle (2.5pt);
  \fill (1,0) circle (2.5pt);
  \fill (1.5,1) circle (2.5pt);
  \fill (2,2) circle (2.5pt);
  \fill (2.5,1) circle (2.5pt);
  \fill (3,0) circle (2.5pt);
  \fill (3.5,1) circle (2.5pt);
  \fill (4,2) circle (2.5pt);
  \fill (4.5,3) circle (2.5pt);
  \fill (5,2) circle (2.5pt);
  \fill (5.5,1) circle (2.5pt);
  \fill (6,0) circle (2.5pt);
  \fill (6.5,1) circle (2.5pt);
  \fill (7,2) circle (2.5pt);
  \fill (7.5,3) circle (2.5pt);

  \node[right] at (8,3) {$h$};
  \node[right] at (8,2) {$2h$};
  \node[right] at (8,1) {$4h$};
  \node[right] at (8,0) {$8h$};
\end{tikzpicture}
\caption{$\algFMG$}
\end{subfigure}

\caption{Schematic of three basic multigrid cycles}
\label{fig:mgcycles}
\end{figure}

The {\em basic} V- and $W$-cycling schemes are depicted schematically in Figure~\ref{fig:mgcycles} (a) and (b), respectively. The FMG-cycling scheme depicted in Figure~\ref{fig:mgcycles} (c) is described in Section~\ref{sec:FMG}. Other cycling schemes such as those introduced in \cite{Avnat2023} exist, but are used less frequently in practice. Before we discuss FMG as a {\em direct} solver in the PDE sense, we first discuss the complexity of the basic cycles as {\em iterative} solvers.

\section{Iterative multigrid solver}
\label{sec:ims}

Accurate error measures are difficult to obtain for most iterative solvers applied to discretized elliptic PDEs. Convergence estimates typically use the evolving residuals or similar quantities that unfortunately hide smooth errors that must be accounted for at an additional cost that grows with the matrix condition number. Single-grid methods applied to large-scale problems typically exacerbate this difficulty because they tend to produce errors that are predominantly smooth. This difficulty is one of the reasons that sharp {\em a posteriori} error results are lacking in many applications. 

An {\em iterative} multigrid solver (IMG) that iterates with one of the cycling schemes discussed thus far is perhaps less troublesome in this sense because it tends to balance the errors better in terms of smooth and oscillatory components. Hidden smooth error is therefore less of an issue. However, it is still somewhat of a concern for IMG because a suitable balance between the errors may be difficult to guarantee. The real difficulty is to devise an effective stopping criterion. This might be expressed by the need to guarantee a certain level of energy-norm accuracy when the residual norm falls below an $\epsilon$ that must somehow be determined. This residual-norm stopping criterion is what the particular IMG given below in Algorithm~\ref{alg:IMG} uses. It is based on a generic multigrid method denoted by $\algMG^h\left({\vec v}^h, {\vec g}^h,\mu,\nu,h\right)$ such as a V-cycle or W-cycle. We apply IMG directly to equation (\ref{equation}) to emphasize its role as a solver, but start the algorithm by setting ${\vec g}^h$ as the residual so that multigrid cycles are nevertheless applied to (\ref{resid-eq}). Since these cycles are not provided with a better initial guess for the correction ${\vec v}^h$ to the current approximation ${\vec u}^h$ to the solution of (\ref{equation}), they start with ${\vec v}^h = {\vec 0}^h$.

This IMG scheme is expressed by
\begin{equation}
{\vec u}^h \leftarrow \algIMG^h\left({\vec u}^h, {\vec g}^h,\mu,\nu,h,\epsilon\right)
\end{equation}
and defined in Algorithm~\ref{alg:IMG}.

\noindent
\begin{center}
\begin{minipage}{1.\linewidth}
\begin{algorithm}[H]
\caption{$\algIMG^h\left({\vec u}^h, {\vec g}^h,\mu,\nu,h,\epsilon \right)$; Iterative multigrid solver}
\label{alg:IMG}
\begin{algorithmic}[1]
      \Require ${\mat L}^h$, ${\mat R}_h^{2h}$, ${\mat P}_{2h}^h$
      \State ${\vec g}^h = {\mat L}^h {\vec u}^h - {\vec f}^h$
      \Comment{Compute the current residual for (\ref{equation}).}
      \State ${\vec v}^h \leftarrow {\vec 0}^h$
      \Comment{Set initial guess for correction ${\vec v}^h$ to ${\vec 0}^h$.}
      \If{$\|{\mat L}^h {\vec v}^h - {\vec g}^h\|>\epsilon$} 
      \Comment{Iterate on correction ${\vec v}^h$ to ${\vec u}^h$.}
      \State ${\vec v}^h \leftarrow \algMG^h\left({\vec v}^h, {\vec g}^h,\mu,\nu,h\right)$
      \Comment{Apply MG cycle to (\ref{resid-eq}) with current iterate.}
      \EndIf
      \State ${\vec u}^h \leftarrow {\vec u}^h - {\vec v}^h$
      \Comment{Correct approximation ${\vec u}^h$ to solution of (\ref{equation}).}
      \State\Return ${\vec u}^h$
\end{algorithmic}
\end{algorithm}
\end{minipage}
\end{center}
\vspace{4mm}

The aim of IMG based on multiple multigrid cycles like the V-cycle and W-cycle is to achieve optimality by converging uniformly well (that is, with a {\em per-cycle} convergence factor bound that is independent of $h$) at a cost equivalent to just a few relaxation sweeps on the finest level. This optimality is often expressed  {\em theoretically} as $O(n)$, where $n$ is the number of grid points on the finest level. However, this does not mean that these iterations achieve acceptable results at such a cost because the required accuracy increases with $n$ (i.e., as $h$ decreases). That is after all why one would want to have a larger $n$. In other words, the goal of refinement is to increase the accuracy of the discrete approximation to the solution of the PDE, so the number of multigrid cycles must increase accordingly. This typically leads to a {\em total} theoretical cost of at least $O(n\log n)$. 

To see this for the elliptic problems we have in mind, we note first that the discrete solution accuracy typically is $O(h^k)$ for some positive integer k. The solver should attempt to match this order of approximation, which means that refining from grid $h$ to grid $\frac{h}{2}$ should reduce the error by a factor of at least $2^{-k}$. Thus, with a fixed per-cycle energy convergence bound of, say, $\rho <1$, we need at least $\gamma$ cycles to ensure that $\rho^\gamma \le 2^{-k}$, that is, $\gamma \ge \frac{k}{|\log_2 \rho|} > 0$. So at least one extra cycle is needed after refinement by a factor of 2. In general, this adds a factor of $\log n$ to the number of cycles needed, which confirms that the overall cycling cost is at least $O(n\log n)$.

\begin{rem} {\bf A word about theoretical optimality.} The statement of $O(n)$ per-cycle optimality is only correct in theory, that is, in infinite precision. In practice, finite precision must also be taken into account. Since precision must increase with the higher demands on accuracy that come with increasingly finer levels, so too must the cost of the finer-level operations in a V-cycle. It is therefore more accurate to interpret optimality as the cost equivalent to just a few relaxation sweeps on the finest level. For simplicity, we nevertheless use $O(n)$ here, but it should be kept in mind that this is a theoretical statement. It is also important to keep in mind that IMG optimality is a per-cycle statement, not a claim about overall cost. Theoretical $O(n)$ optimality is, however, achievable by FMG as described next. \end{rem}

\section{FMG}
\label{sec:FMG}

This extra cost of iterative multigrid methods comes from starting the cycles with a naive initial guess. With the iteration starting its cycles on a given level, what better initial guess is there than zero? The answer is that the cycles could start by computing approximations on the coarsest level up through the hierarchy so that, by the time the finest level is reached, a fairly accurate approximation to the solution of (\ref{equation}) has already been computed, and this would be done at minimal cost because of the reduced coarse-grid complexities. An important point here is that, on the way up to the finest level, cycling would be performed on the coarse levels as opposed to relaxation alone. This idea is the basis for FMG, which often obtains accurate results at a truly optimal $O(n)$ total theoretical cost. 

\subsection{Algorithm}Returning to the original problem of solving (\ref{equation}), the FMG algorithm is based on $q$ applications per level of a generic multigrid cycle denoted by $\algMG^h\left({\vec v}^h, {\vec g}^h,\mu,\nu,h\right)$, which could either be $\algMV$, $\algMW,$ or any other basic choice. We denote the FMG algorithm by the expression
\begin{equation}
{\vec u}^h \leftarrow \algFMG^h\left(\mu,\nu,q,h\right)
\end{equation}
and define it recursively in Algorithm~\ref{alg:FMG}.

\noindent
\begin{center}
\begin{minipage}{1.\linewidth}
\begin{algorithm}[H]
\caption{$\algFMG^h\left(\mu,\nu,q,h\right)$; Full multigrid}
\label{alg:FMG}
\begin{algorithmic}[1]
      \Require ${\mat P}_{2h}^h$, ${\vec f}^h$
      \If{$h<H$}
      \State ${\vec u}^{2h} \leftarrow \algFMG^{2h}\left(\mu,\nu,q,2h\right)$
      \COMMENT{Apply $\algFMG^{2h}$.}
      \State ${\vec u}^h \leftarrow {\mat P}_{2h}^h{\vec u}^{2h}$
      \COMMENT{Start $\algFMG^h$ with grid $2h$ result.}
      \Else
      \State ${\vec u}^h \leftarrow {\vec 0}^h$
      \COMMENT{Set initial grid $H$ guess to zero.}
      \EndIf
        \State $i \leftarrow 0$
  \While{$i<q$}
       \State ${\vec u}^h \leftarrow \algMG^h\left({\vec u}^h, {\vec f}^h,\mu,\nu,h\right)$
      \COMMENT{Apply $\algMG^h$ $q$ times.}\label{MG}
  \State $i \leftarrow i+1$  
  \EndWhile
       \State\Return ${\vec u}^h$
\end{algorithmic}
\end{algorithm}
\end{minipage}
\end{center}
\vspace{4mm}
See Figure~\ref{fig:mgcycles} (c) for a schematic of FMG using a single V-cycle on each level.

The main goal of FMG is to obtain discretization-level accuracy by using coarse grids in a multilevel hierarchy to provide a good enough initial guess that only a few multigrid cycles are required on the finest grid. More precisely, it attempts to achieve an {\em upper bound} on the error of the final approximation of the exact discrete solution in energy that is of the same order as the {\em upper bound} on the energy error in that exact discrete solution as an approximation to the exact PDE solution. The FMG processes for achieving this objective can be determined by theoretical analyses such as that developed in \cite{mmb}, with a typical theoretically optimal $O(n)$ cost. But a key word here is {\em bound}. FMG does not necessarily reach an approximation to the exact discrete solution that is comparable to the discretization error simply because that discrete solution may be much more accurate than analysis predicts. For example, nothing prevents a discrete solution from being the exact PDE solution, and one would not expect FMG to compute the discrete solution exactly. However, such a case is perhaps not so common in practice, and accuracy at the level of the estimated discretization error might be all that one can expect in any case. 

\subsection{Capabilities}As we suggested in section~\ref{sec:ims}, an important benefit of FMG that is often overlooked is its ability to avoid the pitfalls of practical error estimation. While iterative multigrid solvers tend to produce better estimates of the errors than most other solvers, FMG is in an even better position: when properly designed, it can deliver approximations with errors that are comparable to discretization-error bounds without appealing to any convergence estimate at all. It might be argued that solvers should provide error measures, and that is certainly an important goal. However, one can think of FMG as being integral to the discretization in that it aims to achieve the bound that the discretization promises. Since discretization-error trends are usually determined a priori and are not always estimated in the computation itself, it would be in concert to have FMG and discretization integrated in the overall approach. If, however, an error measure is desired at execution, then the errors due to the solver and the discretization can be estimated according to the coarse-grid process described next.

While FMG tends to rely on bounds as opposed to actual errors, discrete solutions often behave in a consistent pattern with errors that are roughly proportional to a power of the mesh size, that is, $Ch^k$, where $C$ is an unknown constant while $k$ is typically a known integer. This may not necessarily occur on coarse levels, but in many cases it shows up asymptotically with decreasing mesh size, what has been called {\em saturation}. The hierarchy of grids that FMG uses provides an opportunity for estimating such occurrences so that discretization-level accuracy can be obtained. The assumption here is that coarse-grid computation is very inexpensive compared to the cost of multigrid cycles on the finest levels, so substantially more effort can be applied at negligible cost there. (Such is not  the case in a parallel processing mode, however, when the number of processors is on the same order as the number of finest-level unknowns.) 

More specifically, three coarse levels can be subjected to enough cycles to ensure that the resulting approximations are very close to the corresponding discrete solutions. This would also provide estimates of the energy per-cycle convergence factors. The two coarser level results could then be compared in the energy norm to the finer level result to estimate the constant $C$ in the error $Ch^k$ if $k$ is known. Otherwise, four levels could be used to estimate both $C$ and $k$ on very coarse levels, yielding estimates for the energy errors on finer levels. Together with estimates for the convergence factors and arithmetic costs, this would provide a basis for deciding whether to apply more multigrid cycles on a given level, to proceed to the next finer grid, or terminate the process. For example, if the error estimates suggest that the accuracy is not yet acceptable on the finest level, then the choice would be to either work harder there or proceed to the next finer level. That choice could be made on what the estimates suggest about an estimated cost-benefit between the two options. 

In other words, instead of relying on error bounds that may or may not be sharp or even available, relatively inexpensive coarse-grid cycles can be used to decide if adequate accuracy has been achieved or more work is required on the current or next finer grid. This refinement strategy does not of course guarantee that the process is fully able to accurately approximate the actual solutions on finer levels, but it may offer better results in many circumstances. This strategy might be further enhanced by analysis beforehand of the data (e.g., coefficients, source terms, boundary conditions, and domain) to estimate when saturation might occur.

We conclude with two FMG applications that are not always well understood. While the discussions do not go into much any detail, they hopefully suggest additional advantages that FMG may possess. \\

{\bf Eigenproblems.}
In terms of saturation, it should be noted that large-scale eigenproblems have historically been considered to be somewhat harder to solve than matrix equations. After all, eigenvalues and their vectors are typically computed by solving many matrix equations (e.g., by inverse iteration). But this is not so with multigrid for the minimal eigenvalues of elliptic PDEs. In fact, in theory, they are generally easier than matrix equations to solve by multilevel methods, in part because the solutions are as geometrically smooth as possible. But there is more to it than that because the accuracy on the levels generally behaves very regularly in the sense that it is more or less equal to some constant times a power of $h$, so that sharp error bounds might be more easily obtained. This makes FMG performance very predictable.\\

{\bf Nonlinear problems.}
Another benefit of FMG that is often overlooked in practice is its effectiveness for many nonlinear problems. The key point here is that the solution to the grid $2h$ problem is often well inside the region of attraction of linearization methods, so FMG works virtually as well in these cases as it does for linear problems. With FAS \cite{brandt77} and other multigrid ways to treat nonlinearity, it may be possible to avoid linearization entirely. But full multigrid often makes that concern secondary in the sense that the various ways to treat the nonlinearity often perform equally well when coupled with FMG. This of course assumes that the coarse levels can be efficiently formed and processed in practice, which is admittedly an issue with multigrid that must be carefully addressed.

\section{Conclusion}
\label{sec:summary}

At this point, we have introduced the basic ideas behind multigrid and its various cycling schemes. In summary, there are four key principles underlying multigrid: minimization, relaxation, coarsening, and the use of a hierarchy. These principles are helpful to keep in mind in the use and construction of both basic as well as more complicated multigrid algorithms. Specifically:

\begin{itemize}
\item {\bf Minimization}. For symmetric positive definite matrix equations, the energy minimization principle has several distinct advantages that arise from its equivalence to the matrix equation itself. Often, this principle is how the original problem was originally formed. It also controls the Euclidean norm, at least for PDEs that we have in mind, because the energy norm applied to admissible PDE functions is equivalent to the norm that includes the function and its individual partial derivatives. Moreover, as opposed to the Euclidean norm, energy minimization also controls artificial oscillations in the approximations.
\item {\bf Relaxation}. Relaxation used on all levels of the grid hierarchy should be relatively inexpensive so that the cost of the multigrid method is controlled, preferably with numerical complexity on the order of the number of degrees of freedom on the finest grid. It is also crucial to the success of multigrid that the errors that are slow to converge under relaxation are understood well enough to support coarsening. For linear iterative relaxation methods, this means that near-kernel components of the error propagation matrix hopefully exhibit an identifiable local pattern of the errors. For Krylov relaxation in particular, algebraic smoothness means small residuals and that can imply that the errors are nearly constant in a neighborhood defined by strong matrix connections.
\item {\bf Coarsening}. Coarsening should target algebraically smooth errors so that it can work in tandem with relaxation to effectively reduce all errors. Algebraic smoothness implies geometric smoothness in some cases, which means that the fewer points of a $2h$ grid are enough to allow adequate approximation by standard interpolation of the errors left by a few relaxation sweeps. In general, a full understanding of algebraically smooth errors is necessary to develop ways to approximate them on a substantially coarsened level of the fine grid.
\item {\bf Hierarchy}. The presence of a hierarchy of equations from the very coarsest to the finest grid offers advantages beyond the promise of a fast solver.  For one, the hierarchy enables nested iteration that uses coarse levels initially to supply a good initial guess to the solution on the finest grid. Another advantage is the possibility of FMG to deal with nonlinearity effectively. Other advantages stem from the ability to aggressively assess numerical performance on relatively inexpensive coarse levels as a way to determine the overall FMG processes. In any case, nesting multigrid cycles allows FMG to deliver acceptable accuracy with a cost equivalent of just a few finest-grid relaxation sweeps. \\
\end{itemize}

Beyond the ideas presented here, a more advanced Primer \cite{mccormick26b} develops the underlying theory, with principles that suggest how to extend multigrid to more complicated problems, including algebraic multigrid (AMG).

\section*{Acknowledgement}

The authors would like to thank Rob Falgout for his effort in making Figure~\ref{fig:fig}.

\bibliographystyle{plainurl}
\bibliography{mgprimer}

\begin{thebibliography}{10}

\bibitem{Avnat2023}
Or~Avnat and Irad Yavneh.
\newblock {On the Recursive Structure of Multigrid Cycles}.
\newblock {\em SIAM Journal on Scientific Computing}, 45(3):S103--S126, 2023.
\newblock \href {https://doi.org/10.1137/21M1433502}
  {\path{doi:10.1137/21M1433502}}.

\bibitem{brandt77}
Achi Brandt.
\newblock Multi-level adaptive solutions to boundary-value problems.
\newblock {\em Mathematics of Computation}, 31(138):333--390, 1977.
\newblock \href {https://doi.org/10.1090/S0025-5718-1977-0431719-X}
  {\path{doi:10.1090/S0025-5718-1977-0431719-X}}.

\bibitem{briggs00}
W.~Briggs, V.~Henson, and S.~McCormick.
\newblock {\em {A Multigrid Tutorial, Second Edition}}.
\newblock Society for Industrial and Applied Mathematics, 2000.
\newblock \href {https://doi.org/10.1137/1.9780898719505}
  {\path{doi:10.1137/1.9780898719505}}.

\bibitem{Golub83}
Gene~H. Golub and Charles F.~Van Loan.
\newblock {\em {Matrix Computations}}.
\newblock The Johns Hopkins University Press, 3rd edition, 1983.

\bibitem{haase-langer}
G.~Haase and U.~Langer.
\newblock {Multigrid methods: from geometrical to algebraic versions}.
\newblock {\em Bourlioux, A., Gander, M.J., Sabidussi, G. (eds) Modern Methods
  in Scientific Computing and Applications, NATO Science Series, Springer}, 75,
  2002.
\newblock \href {https://doi.org/10.1007/978-94-010-0510-4_4}
  {\path{doi:10.1007/978-94-010-0510-4_4}}.

\bibitem{mmb}
J.~Mandel, S.~McCormick, and R.~Bank.
\newblock Variational multigrid theory.
\newblock In Stephen~F. McCormick, editor, {\em Multigrid Methods}, chapter~5,
  pages 131--177. Society for Industrial and Applied Mathematics, 1987.
\newblock \href {https://doi.org/10.1137/1.9781611971057.ch5}
  {\path{doi:10.1137/1.9781611971057.ch5}}.

\bibitem{mccormick84}
Stephen McCormick.
\newblock {Multigrid Methods for Variational Problems: Further Results}.
\newblock {\em SIAM Journal on Numerical Analysis}, 21(2):255--263, 1984.
\newblock \href {https://doi.org/10.1137/0721018} {\path{doi:10.1137/0721018}}.

\bibitem{mccormick26b}
Stephen~F. McCormick and Rasmus Tamstorf.
\newblock {Multigrid Primer: Theory, Principles, and Algebraic Multigrid}.
\newblock Technical report, The University of Colorado at Boulder, CU-Boulder,
  2026.
\newblock URL: \url{https://amath.colorado.edu/~stevem/mgprimer_amg.pdf}.

\bibitem{strang-fix}
G.~Strang and G.~J. Fix.
\newblock {\em An Analysis of the Finite Element Method}.
\newblock Prentice-Hall, Englewood Cliffs, NJ, 1973.

\bibitem{Tamstorf2015}
Rasmus Tamstorf, Toby Jones, and Stephen~F. McCormick.
\newblock Smoothed aggregation multigrid for cloth simulation.
\newblock {\em ACM Transactions on Graphics}, 34(6), Oct 2015.
\newblock \href {https://doi.org/10.1145/2816795.2818081}
  {\path{doi:10.1145/2816795.2818081}}.

\bibitem{trottenberg00}
Ulrich Trottenberg, Cornelius~W. Oosterlee, and Anton Schuller.
\newblock {\em Multigrid}.
\newblock Academic Press, 2000.

\bibitem{Vassilevski2008}
Panayot~S. Vassilevski.
\newblock {\em Multilevel Block Factorization Preconditioners: Matrix-based
  Analysis and Algorithms for Solving Finite Element Equations}, volume~68 of
  {\em Lecture Notes in Computational Science and Engineering}.
\newblock Springer, New York, 2008.
\newblock \href {https://doi.org/10.1007/978-0-387-71564-3}
  {\path{doi:10.1007/978-0-387-71564-3}}.

\bibitem{vassilevski2010}
Panayot~S. Vassilevski.
\newblock Approximation properties of coarse spaces by algebraic multigrid.
\newblock Technical Report LLNL-PROC-422574, Lawrence Livermore National
  Labratory, 2010.
\newblock URL: \url{http://www.osti.gov/servlets/purl/1114745/}.

\bibitem{xu2017}
J.~Xu and L.~Zikatanov.
\newblock Acceleration of convergence of a two-level algorithm by smoothing
  transfer operator.
\newblock {\em Acta Numerica}, pages 1--127, 2017.
\newblock \href {https://doi.org/10.1017/S09624929}
  {\path{doi:10.1017/S09624929}}.

\end{thebibliography}

\end{document}